\input amstex


\def\b1{\text{\bf 1}}

\def\BC{{\Bbb C}}

\def\BP{{\Bbb P}}
\def\BZ{{\Bbb Z}}
\def\BR{{\Bbb R}}
\def\BT{{\Bbb T}}
\def\BW{{\Bbb W}}
\def\CA{{\Cal A}}

\def\CL{{\Cal L}}

\def\CT{{\Cal T}}

\def\gth{{\frak h}}



\def\iso{\buildrel\sim\over\longrightarrow}

\parskip=6pt

\documentstyle{amsppt}
\document
\NoBlackBoxes


\centerline{\bf  Deformations of chiral algebras} 

\centerline{\bf and quantum cohomology of toric varieties}

\bigskip

\centerline{Fyodor Malikov\footnote{partially supported by an
NSF grant} and Vadim Schechtman}

\bigskip\bigskip

Let $X$ be a smooth complex variety. 
It was shown in [MSV]  that the complex cohomology algebra  
$H^*(X)$ may be obtained as a cohomology of a certain vertex  
algebra $H^{ch}(X)$  canonically associated with $X$. 
By definition, $H^{ch}(X)=H^*(X;\Omega^{ch}_X)$, where  $\Omega^{ch}_X$ 
is a sheaf of vertex superalgebras
 constructed in [MSV]. 
(If $X$ is compact, then $H^{ch}(X)$ may be called the {\it chiral 
Hodge cohomology} algebra of $X$.) The algebra $H^{ch}(X)$ is equipped 
with a canonical odd derivation $Q$ of square zero, 
and the cohomology of $H^{ch}(X)$ 
with respect to $Q$ is equal to $H^*(X)$. 

In the very interesting paper [B] Borisov defined for a toric complete 
intersection $X$  a certain vertex superalgebra $V(X)$ equipped with 
an odd derivation of square zero so that $H^{ch}(X)$ equals 
the cohomology of $V(X)$ with respect to this derivation.   
It follows that $H^*(X)$ may also be 
represented as the cohomology of $V(X)$ with respect to another
odd derivation $d$.

Let $X$ be a smooth complete toric variety. 
In the present note we include Borisov's algebra $V(X)$ 
and its derivation $d$ in a family $(V_q(X),d_q)$ 
of vertex superalgebras with derivation, parametrized by 
$q\in H^2(X)$, so that the cohomology of $V_q(X)$ with 
respect to $d_q$ is equal to the {\it quantum cohomology} algebra 
of $X$.  

In sect. 2.5 we present a simpler version of this construction in the
case of $\BP^{N}$ and apply the deformation technique to compute
$H^*(\BP^{N};\Omega^{ch}_{\BP^{N}})$

We also get similar (partial) results for Fano hypersurfaces 
in $P^N$.

 \bigskip\bigskip

\centerline{\bf \S 1. Borisov's construction} 

\bigskip\bigskip

{\bf 1.1.} {\it Lattice vertex algebras.} Let $L$ be a free abelian
group on $2N$ generators $A^{i},B^{i},\;1\leq i\leq N$. Give $L$ an integral
lattice structure by defining a bilinear symmetric $\BZ$-valued form

$$ (.,.): L\times L\rightarrow \BZ$$

so that

$$ (A^{i},B^{j})=\delta_{ij},\;(A^{i},A^{j})=(B^{i},B^{j})=0.$$

Introduce the complexification of $L$: 

$$\gth_{L}=L\otimes_{\BZ}\BC.$$

The bilinear form $(.,.)$ carries over to $\gth_{L}$ by bilinearity. Let

$$ \hat{\gth_{L}}=\gth_{L}\otimes\BC[t,t^{-1}]\oplus\BC K$$

be a Lie algebra with bracket

$$[x\otimes t^{i},y\otimes t^{j}]=i(x,y)\delta_{i+j} K,\;
[x\otimes t^{i}, K]=0.$$

Associated with $L$ there is a group algebra $\BC[L]$ with basis
$e^{\alpha},\alpha\in L,$ and multipliciation

$$e^{\alpha}\cdot e^{\beta}=e^{\alpha+\beta},\; e^{0}=1,\;
 \alpha,\beta\in L.$$

Denote by $S_{\gth_{L}}$ the symmetric algebra of the space
$\gth_{L}\otimes t^{-1}\BC[t^{-1}]$. The space
$ S_{\gth_{L}}\otimes\BC[L]$
carries the well-known vertex algebra structure, see for example [K].
Borisov proposes to enlarge this lattice vertex algebra by fermions
as follows.

We tacitly assumed that $\gth_{L}$ is a purely even vector space:
$\gth_{L}^{(0)}=\gth_{L},\;\gth_{L}^{(1)}=0$. Let $\Pi\gth_{L}$ satisfy the relations
$\Pi\gth_{L}^{(1)}=\gth_{L},\;\Pi\gth_{L}^{(0)}=0$. Thus $\Pi\gth_{L}$ is a purely odd
vector space with basis $\Psi^{i}, \Phi^{i}$ carrying the following odd
bilinear form:

$$
 (.,.): \Pi\gth_{L}\times \Pi\gth_{L}\rightarrow \BC,
$$

$$ 
(\Psi^{i},\Phi^{j})=\delta_{ij},\;(\Psi^{i},\Psi^{j})=
(\Phi^{i},\Phi^{j})=0.
$$

Given all  this, one defines the Clifford algebra, $Cl_{\gth_{L}}$,
to be the    vector superspace

$$
Cl_{\gth_{L}}=\Pi\gth_{L}\otimes\BC[t,t^{-1}]\oplus\BC K',\;
Cl_{\gth_{L}}^{(1)}=\Pi\gth_{L}\otimes\BC[t,t^{-1}],
Cl_{\gth_{L}}^{(0)}=\BC K',
$$

with (super)bracket $[x\otimes t^{i},y\otimes t^{j}]=(x,y)\delta_{i+j} K'$.

Let $\Lambda_{\gth_{L}}$ be the symmetric algebra of the superspace

$$
\oplus_{i=1}^{N}( \Phi^{i}\otimes\BC[t^{-1}]\oplus
\Psi^{i}\otimes t^{-1}\BC[t^{-1}].
$$
(If  we had been allowed to forget about the parity, we would 
have equivalently
defined $\Lambda_{\gth_{L}}$ to be the exterior algebra of the indicated space.)
The space $\Lambda_{\gth_{L}}$ carries the well-known vertex algebra structure,
see for example [K].

Finally let

$$V_{L}=\Lambda_{\gth_{L}}\otimes S_{\gth_{L}}\otimes\BC[L].$$

Being a tensor product of vertex algebras, $V_{L}$ is also a vertex algebra.

\bigskip

{\bf 1.2.} {\it Explicit description of the vertex algebra structure on
$V_{L}$.} To simplify the notation, we identify $\BC[L]$ with the subspace
$1\otimes 1\otimes\BC[L]$. As an $\hat{\gth_{L}}\oplus Cl_{\gth_{L}}$-module,
$V_{L}$ is a direct sum of irreducibles and there is one irreducible
module, $V_{L}(\alpha)$, for each $\alpha\in L$. $V_{L}(\alpha)$ is
freely generated by the supercommutative associative algebra
$S_{\gth_{L}}\otimes\Lambda_{\gth_{L}}$ from the highest weight vector $e^{\alpha}$.
The words ``highest weight vector'' mean that the following relations
hold:
$$
A_{n}^{i}e^{\alpha}=\Psi^{i}_{n}e^{\alpha}=B_{n}^{i}e^{\alpha}=
\Phi^{i}_{n+1}e^{\alpha}=0,\;n\geq 0,
$$
 
$$
Ke^{\alpha}=K'e^{\alpha}=e^{\alpha},\;
  xe^{\alpha}=(x,\alpha)e^{\alpha},\; x\in\gth_{L}.
$$

Thus, $V_{L}(\alpha),\alpha\in L,$ are different as
$\hat{\gth_{L}}\oplus Cl_{\gth_{L}}$-modules, but  isomorphic as
$\hat{\gth_{L}}_{1}\oplus Cl_{\gth_{L}}$-modules, where 
$\hat{\gth_{L}}_{1}\subset\hat{\gth_{L}}$ is the subalgebra linearly spanned
by $x\otimes t^{i},i\neq 0, x\in\gth_{L}$. In fact, the multiplication by
$e^{\beta}$ provides an isomorphism of 
$\hat{\gth_{L}}_{1}\oplus Cl_{\gth_{L}}$-modules:

$$
e^{\beta}: V_{L}(\alpha)\rightarrow V_{L}(\alpha+\beta),\;
x\otimes e^{\alpha}\mapsto x\otimes e^{\alpha+\beta}.
$$

Let us now define the state-field correspondence, that is, attach a field
$x(z)\in\text{End}(V_{L})((z,z^{-1}))$
 to each state $x\in V_{L}$. As has become customary, we shall write
$x_{i}$ for $x\otimes t^{i}$ ($x\in\gth_{L}\text{ or }\Pi\gth_{L}$). We have:
$$
(x_{-n-1}e^{0})(z)=\frac{1}{n!}x(z)^{(n)},\; x\in\gth_{L},
$$
where
$$
x(z)=\sum_{j\in\BZ}x_{j}z^{-j-1}.
$$
In particular, $(x_{-1}e^{0})(z)=x(z)$. 

We continue in the same vein:
$$
(\Phi_{-n}^{i}e^{0})(z)=\frac{1}{n!}\Phi^{i}(z)^{(n)},
$$
where
$$
\Phi^{i}(z)=\sum_{j\in\BZ}\Phi^{i}_{j}z^{-j};
$$
$$
(\Psi^{i}_{-n-1}e^{0})(z)=\frac{1}{n!}\Psi^{i}(z)^{(n)},
$$
where
$$
\Psi^{i}(z)=\sum_{j\in\BZ}\Psi^{i}_{j}z^{-j-1};
$$
$$
e^{\alpha}(z)=e^{\alpha}\cdot\exp{(-\sum_{n<0}\frac{\alpha_{n}}{n}z^{-n})}
\cdot
\exp{(-\sum_{n>0}\frac{\alpha_{n}}{n}z^{-n})}\cdot z^{\alpha_{0}}.
$$
Finally,
$$
x^{(1)}_{-n_{1}}\cdot x^{(2)}_{-n_{2}}\cdots x^{(k)}_{-n_{k}}\cdot e^{\alpha}(z)
=
:x^{(1)}_{-n_{1}}(z) x^{(2)}_{-n_{2}}(z)\cdots x^{(k)}_{-n_{k}}(z)
 e^{\alpha}(z):.
$$

The vertex algebra structure on $V_{L}$ is equivalently described by the
 following family of $n$-th products ($n\in\BZ$):
$$
_{(n)}: V_{L}\otimes V_{L}\rightarrow V_{L},\; x\otimes y\mapsto
x_{(n)}y\buildrel\text{def}\over = (\int x(z)z^{n})(y),
$$
where $\int x(z)z^{n}$ stands for the linear transformation of $V_{L}$
equal to the coefficient of $z^{-n-1}$ in the series $x(z)$.

\bigskip

{\bf 1.3.} {\it Degeneration of $V_{L}$.} Denote by $L_{A}$ the subgroup of
 $L$
generated by $A^{i},\; i=1,...,N$.
 Any smooth toric variety $X$ can be defined via a fan, $\Sigma$,
that is, a collection of ``cones'' lying in $L_{A}$. Borisov uses such
$\Sigma$ to define a certain degeneration, $V_{L}^{\Sigma}$, of the vertex
algebra structure on $V_{L}$. He further shows that the cohomology of
$V_{L}^{\Sigma}$ with respect to a certain differential
$
D^{\Sigma}:\; V_{L}^{\Sigma}\rightarrow V_{L}^{\Sigma}
$
equals $H^{*}(X,\Omega^{ch}_{X})$, where $\Omega^{ch}_{X}$ is the chiral
de Rham complex of [MSV]. Let us describe the outcome of this construction
in the case when $X=\BP^{N}$.

Consider the following set of $N+1$ elements of $L_{A}$:
$\xi_{1}=A^{1},\xi_{2}=A^{2},...,\xi_{N}=A^{N},
\xi_{N+1}=-A^{1}-A^{2}-\cdots
-A^{N}$. Define the   cone $\Delta_{i}\subset L_{A}$ to be
 the  set of all non-negative integral linear combinations of the elements
$\xi_{1},...,\xi_{i-1},\xi_{i+1},...,\xi_{N+1}$. It is easy to see that
$L_{A}=\cup_{i}\Delta_{i}$ and the intersection $\Delta_{i}\cap\Delta_{j}$ is
a face of both $\Delta_{i}$ and $\Delta_{j}$. The fan $\Sigma$ in this case
is the set consisting of $\Delta_{1},...,\Delta_{N+1}$ and their faces.

We now define $V_{L}^{\Sigma}$ to be a vertex algebra equal to
$V_{L}$ as a vector space with $n$-th product $_{(n),\Sigma}$ as follows:

if $\{\sum_{i}n_{i}A^{i},\sum_{i}n'_{i}A^{i}\}\subset\Delta_{j}$
   for some $j$, then

$$
(x\otimes e^{\sum_{i}m_{i}B^{i}+\sum_{i}n_{i}A^{i}})
_{(n),\Sigma}
(y\otimes e^{\sum_{i}m'_{i}B^{i}+\sum_{i}n'_{i}A^{i}})
$$
$$
= (x\otimes e^{\sum_{i}m_{i}B^{i}+\sum_{i}n_{i}A^{i}})
_{(n)}
(y\otimes e^{\sum_{i}m'_{i}B^{i}+\sum_{i}n'_{i}A^{i}});
$$
   
otherwise
$$
(x\otimes e^{\sum_{i}m_{i}B^{i}+\sum_{i}n_{i}A^{i}})
_{(n),\Sigma}
(y\otimes e^{\sum_{i}m'_{i}B^{i}+\sum_{i}n'_{i}A^{i}})=0,
$$
 
where $_{(n)}$ stands for the n-th product on $V_{L}$.
The fact that these new operations satisfy the Borcherds identities
 can be proved
by including both $V_{L}$ and $V_{L}^{\Sigma}$ in a 1-parameter family
of vertex algebras; this will be done in 2.1.

Let
$$
D=\int\sum_{i=1}^{N}\Psi^{i}(z)(e^{A^{i}}-e^{-\sum_{j}A^{j}})(z).
\eqno{(1.1)}
$$
It is obvious that $D\in\text{End}(V_{L}^{\Sigma})$ and $D^{2}=0$;
 therefore,
the cohomology $H_{D}(V_{L}^{\Sigma})$ arises.

{\bf Theorem 1.3.} ([B])
$$ 
H_{D}(V_{L}^{\Sigma})=H^{*}(\BP^{N},\Omega^{ch}_{\BP^{N}}).
$$

\bigskip\bigskip

\centerline{\bf \S 2. Deforming $H^{*}(\BP^{N})$}

\bigskip\bigskip

{\bf 2.1.} {\it The family $V_{L,q}$.}

Here we exhibit a family of vertex algebras, $V_{L,q}, q\in\BC$, so that
$V_{L,q}$ is isomorphic to $V_{L}$
if $q\neq 0$ and $V_{L,0}$ is isomorphic to $V_{L}^{\Sigma}$; cf.
 the end of sect.8 [B].

Define the height function
$$
ht:\; L_{A}\rightarrow\BZ_{>}
$$ 
as follows. It is easy to see that each $\alpha\in L_{A}$ is uniquely
represented in the form
$$
\alpha=\sum_{i=1}^{N+1}n_{i}\xi_{i}
\eqno{(2.1)}
$$
so that all $n_{i}\geq 0$ and $\#\{i: n_{i}>0\}\leq N$. Let
$$
ht(\alpha)=\sum_{i}n_{i},
$$
where $n_{1},...,n_{N}$ are as in (2.1).

Define the linear automorphism
$$
t_{q}: V_{L}\rightarrow V_{L},\; q\in\BC-\{0\}
$$
by the formula
$$
t_{q}(x\otimes e^{\sum_{i}m_{i}B^{i}+\sum_{i}n_{i}A^{i}} )
=q^{ht(\sum_{i}n_{i}A^{i} )}
x\otimes e^{\sum_{i}m_{i}B^{i}+\sum_{i}n_{i}A^{i}}.
$$
Define $V_{L,q}$ to be the vertex algebra equal to $V_{L}$ as a vector
space with the following n-th product:
$$
(x\otimes e^{\sum_{i}m_{i}B^{i}+\sum_{i}n_{i}A^{i}})_{(n),q}
(y\otimes e^{\sum_{i}m'_{i}B^{i}+\sum_{i}n'_{i}A^{i}})
$$
$$=
t_{q}^{-1}(t_{q}(x\otimes e^{\sum_{i}m_{i}B^{i}+\sum_{i}n_{i}A^{i}} )
_{(n)}t_{q}(y\otimes e^{\sum_{i}m'_{i}B^{i}+\sum_{i}n'_{i}A^{i}} )).
$$
By definition,
$$
t_{q}: V_{L,q}\rightarrow V_{L},\; q\in\BC-\{0\},
$$
is a vertex algebra isomorphism. It is also easy to see that
if  $\sum_{i}n_{i}A^{i}$ and $\sum_{i}n'_{i}A^{i}$ belong to
the same cone from $\Sigma$, then
$$
(x\otimes e^{\sum_{i}m_{i}B^{i}+\sum_{i}n_{i}A^{i}})_{(n),q}
(y\otimes e^{\sum_{i}m'_{i}B^{i}+\sum_{i}n'_{i}A^{i}})
$$
$$=
(x\otimes e^{\sum_{i}m_{i}B^{i}+\sum_{i}n_{i}A^{i}} )
_{(n)}(y\otimes e^{\sum_{i}m'_{i}B^{i}+\sum_{i}n'_{i}A^{i}} );
$$
 otherwise
$$
(x\otimes e^{\sum_{i}m_{i}B^{i}+\sum_{i}n_{i}A^{i}})_{(n),q}
(y\otimes e^{\sum_{i}m'_{i}B^{i}+\sum_{i}n'_{i}A^{i}})
$$
$$\in q\BC[q]
(x\otimes e^{\sum_{i}m_{i}B^{i}+\sum_{i}n_{i}A^{i}} )
_{(n)}(y\otimes e^{\sum_{i}m'_{i}B^{i}+\sum_{i}n'_{i}A^{i}} ).
$$
 Two things follow at once:
first,  the operations 
$$
_{(n),0}=\lim_{q\rightarrow 0}{_{(n),q}}, n\in\BZ
$$ 
are well defined
and satisfy the Borcherds identities; second,  the vertex algebra,
$V_{L,0}$, obtained in this way is isomorphic to $V_{L}^{\Sigma}$. 
By the way,
 this remark proves that $V_{L}^{\Sigma}$ is indeed a vertex algebra.

To get a better feel for this kind of deformation, and for the future use, 
let us
consider the subspace $\BC[L_{A}]\subset V_{L,q}$ with basis $e^{\alpha},
\alpha\in L_{A}$. The $(-1)$-st 
product makes this space  a commutative algebra.
The subspace $\BC[\Delta_{j}]$ defined to be the linear span of
$e^{\alpha},
\alpha\in \Delta_{j},$ is a polynomial ring on generators 
$e^{\xi_{1}},...,e^{\xi_{j-1}},e^{\xi_{j+1}},...,e^{\xi_{N+1}}$.
For example, if we denote $x_{i}=e^{A^{i}}$, then
$\BC[\Delta_{N+1}]=\BC[x_{1},...,x_{N}]$ and this isomorphism identifies
$e^{\sum_{j}n_{j}A^{j}}$ with the monomial
 $x_{1}^{n_{1}}\cdots x_{N}^{n_{N}}$.

The entire  $\BC[L_{A}]$ is not a polynomial ring. For example, as follows
from the definition of the deformation, there is
a relation
$$
(e^{-A^{1}-\cdots -A^{N}})_{(-1)}(e^{A^{1}+\cdots +A^{N}})=q^{N+1} e^{0},
$$
because
 $ht(0)=0$, $ht(A^{1}+\cdots +A^{N})=N$, $ht(-A^{1}-\cdots -A^{N})=1$.
If we let $T=e^{-A^{1}-\cdots -A^{N}}$, then the last equality rewrites
as follows:
$$
Tx_{1}x_{2}\cdots x_{N}=q^{N+1},
$$
and a moment's thought shows that in fact
$$
\BC[L_{A}]=\BC[x_{1},...,x_{N},T]/(Tx_{1}x_{2}\cdots x_{N}-q^{N+1}).
$$

Being a group algebra, $\BC[L_{A}]$ carries another algebra structure, a 
priori
different from the one we just described and independent of $q$.
We see that   the two structures are isomorphic if $q\neq 0$;
  at $q=0$, however, the   one we just described
 degenerates in an algebra with zero divizors.

\bigskip

{\bf 2.2.} {\it The algebra} $H^{*}(\BP^{N})$.

Let
$$
Q(z)=A^{i}(z)\Phi^{i}(z)-\sum_{j}\Phi^{j}(z)',
$$
$$
G(z)=B^{i}(z)\Psi^{i}(z),
$$
$$
J(z)=:\Phi^{i}(z)\Psi^{i}(z):+ \sum_{j}B^{j}(z)',
$$
$$
L(z)=:B^{i}(z)A^{i}(z):+ :\Phi^{i}(z)'\Psi^{i}(z):,
$$
where the summation with respect to repeated indices is assumed.

One checks that the Fourier components of these 4 fields satsify the
commutation relations of the $N=2$ algebra. It is also easy to see that 
the fields
$G(z), L(z)$ commute with Borisov's differential $D$, see (1.1), and therefore
define the fields, to be also denoted $G(z), L(z)$, acting on
$H_{D}(V_{L}^{\Sigma})$.

The fields $Q(z), J(z)$ do not commute with $D$, but their Fourier
 components
$Q_{0}=\int Q(z)$ and $J_{0}=\int J(z)$ do:
$$
[Q_{0},D]=[J_{0},D]=0.
$$
Thus we get 2 operators, to be also denoted $Q_{0}, J_{0}$, acting on
$H_{D}(V_{L}^{\Sigma})$. All  this is summarized by saying that
$H_{D}(V_{L}^{\Sigma})$ is a topological vertex  algebra.

A glance at the formulas on p. 17 of [B] shows that the isomorphism
$H_{D}(V_{L}^{\Sigma})=H^{*}(\BP^{N},\Omega^{ch}_{\BP^{N}})$
(see Theorem 1.3) identifies these $G(z),L(z),Q_{0},J_{0}$ with the fields
(operators) constructed in [MSV] and denoted in the same way.
 One of the main results of [MSV] then gives
$$
H^{*}(\BP^{N})=H_{Q_{0}}(H_{D}(V_{L}^{\Sigma})).
\eqno{(2.2)}
$$
Further, the algebra structure of $H^{*}(\BP^{N})$ is restored from the
 $(-1)$-st product on $H_{D}(V_{L}^{\Sigma})$.

\bigskip

{\bf 2.3.} {\it Deformation of the algebra structure.}

It follows from the proof of Theorem 2.3 below
 that the cohomology (2.2) can be calculated in the reversed
order:
$$
H^{*}(\BP^{N})=H_{D}(H_{Q_{0}}(V_{L}^{\Sigma})).
\eqno{(2.3)}
$$
Note that $D$ and $Q_{0}$ can also be regarded as well-defined operators
acting on the deformed algebra:
$$
D=\int\sum_{i=1}^{N}\Psi^{i}(z)(e^{A^{i}}-e^{-\sum_{j}A^{j}})(z),
Q_{0}=\int A^{i}(z)\Phi^{i}(z):\; V_{L,q}\rightarrow V_{L,q}.
$$
It is immediate to see that $D^{2}=Q_{0}^{2}=0$ on $V_{L,q}$ as well.
Moreover,
$$
[D,Q_{0}]=0.
\eqno{(2.4)}
$$
Indeed, the formulas of 1.2 imply the following OPE:
$$
\sum_{i=1}^{N}\Psi^{i}(z)(e^{A^{i}}-e^{-\sum_{j}A^{j}})(z)\cdot
A^{j}(w)\Phi^{j}(w)=
\frac{\sum_{j}e^{A^{j}}(w)'-e^{-\sum_{j}A^{j}}(w)'}{z-w}.
$$
Therefore,
$$
[D,Q_{0}]=\int \{\sum_{j}e^{A^{j}}(w)-e^{-\sum_{j}A^{j}}(w)\}'=0.
$$

Thus it is  natural to take the
space $H_{D}(H_{Q_{0}}(V_{L,q}))$ for a deformation of $H^{*}(\BP^{N})$.

{\bf Theorem 2.3.}
$$
H_{D}(H_{Q_{0}}(V_{L,q}))=\BC[T]/(T^{N+1}-q^{N+1}).
$$

\bigskip

{\bf Proof.}

1) {\it Computation of} $H_{Q_{0}}(V_{L,q})$. By definition
$$
Q_{0}=\sum_{n\in\BZ}A^{i}_{-n}\Phi^{i}_{n}
\eqno{(2.5)}
$$
Therefore,
$$
[Q_{0},\Psi^{j}_{0}]=A^{j}_{0},\; [Q_{0},G_{0}]=L_{0}.
$$
These relations imply that 
$$
H_{Q_{0}}(V_{L,q})=H_{Q_{0}}(\cap_{j}Ker A^{j}_{0}\cap Ker L_{0}).
\eqno{(2.6)}
$$

It follows from 1.2 that
the space $\cap_{j}Ker A^{j}_{0}\cap Ker L_{0}$ is a linear span of elements
of the form:
$$
\Phi^{i_{1}}_{0}\cdots\Phi^{i_{m}}_{0}e^{\sum_{i}n_{i}A^{i}}.
$$
Formula (2.5) shows that the restriction of $Q_{0}$ to this subspace is 0.
Thus
$$
H_{Q_{0}}(\cap_{j}Ker A^{j}_{0}\cap Ker L_{0})=
\cap_{j}Ker A^{j}_{0}\cap Ker L_{0}.
$$
The $(-1)$-st product makes this subspace  a supercommutative algebra.
In the same way as in 2.1 we get an isomorphism
$$
\cap_{j}Ker A_{0}^{j}\cap Ker L_{0}=\BC[x_{1},...,x_{N},T;\Phi_{1},...,
\Phi_{N}]/(Tx_{1}\cdots x_{N}-q^{N+1}),
$$
where $\Phi_{1},...,
\Phi_{N}$ are understood as grassman variables
 ($[x_{i},\Phi_{j}]=[T,\Phi_{j}]=0$,
$\Phi_{i}\Phi_{j}+\Phi_{j}\Phi_{i}=0$) and $(Tx_{1}\cdots x_{N}-q^{N+1})$
stands for the ideal generated by $Tx_{1}\cdots x_{N}-q^{N+1}$.

2) {\it Computation of $H_{D}(H_{Q_{0}}(V_{L,q}))$}.
 In view of Step 1), we have to restrict $D$ to
$$
\cap_{j}Ker A^{j}_{0}\cap Ker L_{0}. 
$$
The isomorphism
$$
\cap_{j}Ker A_{0}^{j}\cap Ker L_{0}=\BC[x_{1},...,x_{N},T;\Phi_{1},...,
\Phi_{N}]/(Tx_{1}\cdots x_{N}-q^{N+1}),
$$
identifies $D$ with $\sum_{i}(x_{i}-T)\partial/\partial(\Phi_{i})$.
Therefore, the complex
$$
(\BC[x_{1},...,x_{N},T;\Phi_{1},...,
\Phi_{N}]/(Tx_{1}\cdots x_{N}-q^{N+1}),D)
$$
is simply the Koszul resolution of the algebra
$$
\{\BC[x_{1},...,x_{N},T ]/(Tx_{1}\cdots x_{N}-q^{N+1})\}/
(x_{1}-T,x_{2}-T,...,x_{N}-T)
$$
associated with the sequence $x_{1}-T,x_{2}-T,...,x_{N}-T$. This sequence is
regular and we get at once
$$
H_{D}(H_{Q_{0}}(V_{L,q}))=\BC[T]/(T^{N+1}-q^{N+1}). \qed
$$

It is easy to infer from Borisov's proof of Theorem 1.3 that the element
$T=e^{-A^{1}-\cdots -A^{N}}\in V_{L,q}$ is a cocycle representing the
cohomology class proportional to that of a hyperplane in $\BP^{N}$.
This means that the deformation of $H^{*}(\BP^{N})$ we obtained coincides
with the standard one, except that for some reason $q$ happened to be
 raised to the power of $N$.

\bigskip

{\bf 2.4.} {\it Reduction to a single differential.}

Of course it would be nicer to get $H^{*}(\BP^{N})$, or its
deformation, as the cohomology of this or that vertex algebra with respect
to a single differential rather than to compute a repeated cohomology.

{\bf Theorem 2.4.}
$$
H_{D+Q_{0}}(V_{L,q}))=\BC[T]/(T^{N+1}-q^{N+1}).
$$

It is no wonder, in view of Theorem 2.3, that this assertion is a result
of computation of a certain spectral sequence. We shall use several spectral
sequences arising in the following situation, which  is slightly different from
the standard one. Let 
$$
W=\oplus_{n=-\infty}^{+\infty} W^{n}
$$
be a graded vector space with two commuting differentials
$$
d_{1}:\; W^{n}\rightarrow W^{n+1},
d_{2}:\; W^{n}\rightarrow W^{n-1}.
\eqno{(2.7)}
$$
There arise the total differential $d=d_{1}+d_{2}$ and
 the cohomology $H_{d_{1}+d_{2}}(W)$. Note that this cohomology group is
not graded since $d_{1}$ and $d_{2}$ map in opposite directions.
 We can, however,
introduce the filtration
$$
W=\cup_{n} W^{\leq n},\; W^{\leq n}=\oplus_{m=-\infty}^{n} W^{m}.
$$
Then
$$
(d_{1}+d_{2})(W^{\leq n})\subseteq W^{\leq (n+1)}
$$
and there arises a filtration $H_{d_{1}+d_{2}}(W)^{\leq n}$ on the cohomology
and the graded object $Gr H_{d_{1}+d_{2}}(W)$.

It is straightforward to define a spectral sequence
$$
\{ E(W)^{n}_{r}, d^{(r)}:E(W)^{n}_{r}\rightarrow E(W)^{n-r+1}_{r}\},\;
E(W)^{n}_{r+1}= H_{d^{(r)}}(E(W)^{n}_{r}),
\eqno{(2.8)}
$$
the first three terms being as follows:
$$
 E(W)^{n}_{0}=W^{n},\; E(W)^{n}_{1}=H_{d_{1}}(W^{n}),\;
E(W)^{n}_{2}=H_{d_{2}}(H_{d_{1}}(W^{n})),
\eqno{(2.9)}
$$
where
$$
H_{d_{1}}(W^{n})=\frac{\text{Ker}\{d_{1}: W^{n}\rightarrow W^{n+1}\}}
{\text{Im}\{d_{1}: W^{n-1}\rightarrow W^{n}\}},
$$
$$
H_{d_{2}}(H_{d_{1}}(W^{n}))
=\frac{\text{Ker}\{d_{2}: H_{d_{1}}(W^{n})\rightarrow H_{d_{1}}(W^{n-1}) \}}
{\text{Im}\{d_{2}: H_{d_{1}}(W^{n+1})\rightarrow H_{d_{1}}(W^{n})\}}.
$$

In the situation pertaining Theorem 2.4 we take $V_{L,q}$ for $W$,
$Q_{0}$ for $d_{1}$, and $D$ for $d_{2}$. The space $V_{L,q}$ is graded by
fermionic charge; this grading is defined by letting the 
degree of $\Psi^{i}_{j}$
be equal $-1$, the degree of $\Phi^{i}_{j}$ be equal $1$, the degree of
$A^{i}_{j}, B^{i}_{j}, e^{\alpha}$ be equal $0$. By definition, 
$$
Q_{0}(V_{L,q}^{n})\subseteq V_{L,q}^{n+1},
D(V_{L,q}^{n})\subseteq V_{L,q}^{n-1}
$$
and we get a spectral sequence $\{ E(V_{L,q})^{n}_{r}, d^{(r)}\}$.
 
Observe that
the grading by fermionic charge and the corresponding filtration
are infinite in both directions. Therefore, the standard finiteness
conditions that guarantee convergence of spectral sequences fail.
Nevertheless the following lemma holds true.

{\bf Lemma 2.4.} The spectral sequence $\{ E(V_{L,q})^{n}_{r}, d^{(r)}\}$
converges to $H_{Q_{0}+D}(V_{L,q})$ and collapses:
$$
 H_{D}(H_{Q_{0}}(V_{L,q}))=H_{Q+D}(V_{L,q}).
$$

Lemma 2.4 combined with Theorem 2.3 gives Theorem 2.4 at once and
 it remains to prove Lemma 2.4.

{\it Proof of Lemma 2.4.}   Introduce yet another grading of the space
$V_{L,q}^{n}$ as follows. Let $\alpha=(\alpha_{1},...,\alpha_{N+1})$
be an element of the group $\BZ^{N+1}$. Let
$$
V_{L,q}^{n}[\alpha]= (\cap_{i=1}^{N}Ker(A^{i}_{0}-\alpha_{i} Id))\cap
Ker(L_{0}-\alpha_{N+1} Id).
$$
Of course
$$
V_{L,q}^{n}=\oplus_{\alpha\in\BZ^{N+1}}V_{L,q}^{n}[\alpha]
$$
and both the differentials preserve this grading. Therefore all  calculations
can be carried out inside $V_{L,q}^{n}[\alpha]$ with a fixed $\alpha$. Consider
the following two cases.

1) $\alpha\neq 0$. In this case, as was observed in the beginnning of the proof
of Theorem 2.3 (see e.g. (2.6)),
 $H_{Q_{0}}(V_{L,q}[\alpha])=0$ and, therefore,
$E[\alpha]_{1}=0$. It remains to show that $H_{Q_{0}+D}(V_{L,q}[\alpha])=0$.
Let $x\in V_{L,q}[\alpha]^{\leq n}$ be a cocycle. This means that there is
 a ``chain''
of elements $x_{i}\in V_{L,q}[\alpha]^{n-2i} , i=0,1,...,k$ so that
$$
x=\sum_{i=0}^{k}x_{i},
$$
 and the following holds
$$
Q_{0}(x_{0})=0, Q_{0}(x_{i+1})+D(x_{i})=0, D(x_{k})=0, i=0,...,k-1.
\eqno{(2.10)}
$$

We now repeatedly use the condition $H_{Q_{0}}(V_{L,q}[\alpha])=0$ and  
 (2.10) to construct another chain 
$y_{i}\in V_{L,q}[\alpha]^{n-2i-1}, i\geq 0$,
satisfying
$$
Q_{0}(y_{0})=x_{0}, Q_{0}(y_{i+1})+D(y_{i})=x_{i+1}.
\eqno{(2.11)}
$$

Indeed,
since $Q_{0}(x_{0})=0$, there is $y_{0}\in V_{L,q}[\alpha]^{n-1}$ so that
$Q_{0}(y_{0})=x_{0}$.

Since
$$
Q_{0}(-D(y_{0})+x_{1})=DQ_{0}(y_{0})+Q_{0}(x_{1})=D(x_{0})+Q_{0}(x_{1})=0,
$$
 there is $y_{1}\in V_{L,q}[\alpha]^{n-3}$ so that 
$Q_{0}(y_{1})+D(y_{0})=x_{1}$.

In general, having found $y_{i}\in V_{L,q}[\alpha]^{n-2i-1}, 
y_{i-1}\in V_{L,q}[\alpha]^{n-2i+1}$ so that
 $Q_{0}(y_{i})+D(y_{i-1})=x_{i}$,
we calculate as follows:
$$
0=D(0)=D(Q_{0}(y_{i})+D(y_{i-1})-x_{i})=DQ_{0}(y_{i})-D(x_{i}).
$$
Due to (2.10), the last expression rewrites as 
$DQ_{0}(y_{i})+Q_{0}(x_{i+1})$
and we get
$$
-Q_{0}D(y_{i})+Q_{0}(x_{i+1})=0.
$$

Therefore, $Q_{0}(D(y_{i})-x_{i+1})=0$ and there is  $y_{i+1}$ so that
$Q_{0}(y_{i+1})=-D(y_{i})+x_{i+1}$ as desired. 

Formally, (2.11) means that 
$$
(D+Q_{0})(\sum_{i=0}^{\infty}y_{i})=x 
$$
and what does not allow us to conclude immediately that 
$x=\sum_{i=0}^{\infty}x_{i}$ is a coboundary
is that the sum $\sum_{i=0}^{\infty}y_{i}$ looks infinite.   To complete
case 1) it remains to show that
  $y_{i}=0$ for all sufficiently
large $i$. This is achieved by the following dimensional argument.
Note that by construction  
$$
y_{i}\in \oplus_{|m_{j}|<ki}(S_{\gth_{L}}\otimes\Lambda_{\gth_{L}}^{n-2i-1}
\otimes
e^{\sum_{j}m_{j}A^{j}+\sum_{j}\alpha_{j}B^{j}}),
\eqno{(2.12)}
$$
where $k$ is a number independent of $i$. Indeed, each application
of $D$ changes $m_{j}$ by at most 1,  $Q_{0}$ preserves $m_{j}$, and the
linear estimate of $m_{j}$ follows. On the other hand we have an 
explicit formula for $L_{0}$ (see the beginning of 2.2),
 and this formula implies that the smallest eigenvalue
of $L_{0}$ restricted to $\Lambda_{\gth_{L}}^{n-2i-1}$ is nonnegative and
grows faster than a polynomial of degree  2, say $q(i)$,
as $i\rightarrow +\infty$.
The same formula gives
$$
L_{0}e^{\sum_{j}m_{j}A^{j}+\sum_{j}\alpha_{j}B^{j}}=
\sum_{j}m_{j}\alpha_{j} e^{\sum_{j}m_{j}A^{j}+\sum_{j}\alpha_{j}B^{j}}.
$$
Therefore, if $y_{i}\neq 0$, then it is a sum of eigenvectors
associated to eigenvalues of $L_{0}$ greater or equal 
$q(i)-(\alpha_{1}+\cdots\alpha_{n})ki$. Since
this number tends to $+\infty$ as $i\rightarrow +\infty$, we arrive at
contradiction with the assumption $L_{0}y_{i}=\alpha_{N+1}$ if $i$ is
sufficiently large. 
Hence, $y_{i}=0$ for all sufficiently large $i$,
  each cocycle is a coboundary, and
case 1) is accomplished.

2) $\alpha=0$. As we saw in the beginning of the proof of Theorem 2.3,
the restriction of $Q_{0}$ to $V_{L,q}[0]$ is 0 and, by definition, the
complex $(V_{L,q}[0],D+Q_{0})$ is equal to $(E(V_{L,q})[0]_{1}, d^{(1)})$.
$\qed$

\bigskip

{\bf 2.5.} {\it The vertex algebra $H_{D}(V_{L,q})$ and a computation of 
$H^{*}(\BP^{N},\Omega^{ch}_{\BP^{N}})$.}

\bigskip

In this section we prove the following two theorems.

{\bf Theorem 2.5A}  {\it If $q\neq 0$, then $H_{D}(V_{L,q})$ 
 equals the quantum cohomology of $\BP^{N}$.}
\bigskip
{\bf Theorem 2.5B}
{\it
   The natural embedding of sheaves ([MSV], see also  (2.20) below)
$$
\Omega^{*}_{\BP^{N}}\hookrightarrow \Omega^{ch}_{\BP^{N}}
$$
provides an isomorphism 
$$
H^{i}(\BP^{N},\Omega^{*}_{\BP^{N}})\iso H^{i}(\BP^{N},\Omega^{ch}_{\BP^{N}}),\;
0<i<N,
$$
where $\Omega^{*}_{\BP^{N}}$ is the sheaf of all differential forms.}
\bigskip
Recall the  previously known results on the cohomology 
of $\Omega^{ch}_{\BP^{N}}$.
  $\Omega^{ch}_{\BP^{N}}$
is a sheaf of $\widehat{sl}_{N+1}$-modules [MS1], see also  2.5.2. 
In particular, if $U_{0}=\BC^{N}\subset\BP^{N}$ is a big cell, then
$\Gamma(U_{0}, \Omega^{ch}_{\BP^{N}})$ is a generalized Wakimoto module
over $\widehat{sl}_{N+1}$
introduced in [FF]. We proved in [MS1] that
$$
H^{0}(\BP^{N},\Omega^{ch}_{\BP^{N}})=
\Gamma(U_{0}, \Omega^{ch}_{\BP^{N}})^{int},
\eqno{(2.13)}
$$
where $\Gamma(U_{0}, \Omega^{ch}_{\BP^{N}})^{int}$  stands for the
 maximal $sl_{N+1}$-integrable submodule of
$\Gamma(U_{0}, \Omega^{ch}_{\BP^{N}})$.

On the other hand, it follows from the chiral Serre duality [MS2] that
$$
H^{N}(\BP^{N},\Omega^{ch}_{\BP^{N}})=
H^{0}(\BP^{N},\Omega^{ch}_{\BP^{N}})^{d},
\eqno{(2.14)}
$$
where $^{d}$ stands for the
restricted dual.

  Unfortunately, little is known about the structure of 
$\Gamma(U_{0}, \Omega^{ch}_{\BP^{N}})$ and  
$\Gamma(U_{0}, \Omega^{ch}_{\BP^{N}})^{int}$, if $N>1$;
 see, however, [MS1]
for the case of $N=1$.
Otherwise, Theorem 2.5 and (2.13-14) give a   complete
description of $H^{*}(\BP^{N},\Omega^{ch}_{\BP^{N}})$. 

The proofs of  Theorems 2.5A and B are contained in 2.5.2. In  2.5.1
  we collect some well-known material in order
to place these results in the proper context and to formulate (2.18-19), two 
well-known assertions needed in 2.5.2.

\bigskip

{\bf 2.5.1}  
A vertex algebra structure on a vector space $V$ comprises  a countable
 family of   multiplications:
$$
_{(n)}: V\otimes V\rightarrow V,\; x\otimes y\mapsto x_{(n)}y,\; n\in\BZ,
$$
 a map
$$
T: V\rightarrow V,
$$
and a vacuum vector
$$
\b1\in V.
$$
These data satisfy the Borcherds identities which imply, in particular,
that $T$ and $x_{(0)}, x\in V$, are  derivations of
 the $n$-th product for all $n$. Thus,
$$
[T,y_{(j)}]=(Ty)_{(j)},\; [x_{(0)},y_{(j)}]=(x_{(0)}y)_{(j)}.
\eqno{(2.15)}
$$

In the case of the  vertex algebra $V_{L,q}$, the $n$-th multiplication
was defined in the end of 1.2, $\b1$ equals $e^{0}$, and $T$ will be defined
below.

  Call $V$ {\it commutative} (or
holomorphic, see [K] 1.4) if $_{(n)}=0$ for all $n\geq 0$. 
If $V$ is commutative, then the (-1)-st multiplication gives it
the structure of a commutative superalgebra with derivation $T$,
 and the functor arising in this way
is an equivalence of
the category
of commutative vertex algebras and the category of
commutative superalgebras with derivation, see again [K] 1.4. 

If $d_{x}: V\rightarrow V$ is a differential ($d^{2}=0$),
then the cohomology $H_{d_{x}}(V)$ arises. 
We assert that
$$
d_{x}=x_{(0)}\text{ for some }x\in V \Rightarrow 
H_{d_{x}}(V)\text{ is a vertex algebra},
\eqno{(2.16)}
$$
since  all   products on $V$ descend to  $H_{d_{X}}(V)$ due to (2.15).

All vertex algebras we are concerned with are {\it conformal}. This means
that there is a Virasoro field $L(z)=\sum_{i}L_{i}z^{-i-2}$,
$L_{i}\in End (V)$, such that $L_{i}$ satisfy the Virasoro commutation
relations, $T=L_{-1}$,  $L_{0}$ is diagonalizable, and $L(z)$ is the field
attached to the state $L_{-2}\b1\in V$.
The formula at the beginning of 2.2 shows that $V_{L,q}$ is a conformal vertex
algebra, the state $L_{-2}\b1$ being equal to 
$\sum_{i}(B^{i}_{-1}A^{i}_{-1}+\Phi^{i}_{-1}\Psi^{i}_{-1})e^{0}$.
  
 The eigenvalues of 
$L_{0}$ are called conformal weights. Hence a conformal vertex algebra $V$
is graded by conformal weights, $V=\oplus_{n}V_{n}$, and in the case of
$V=V_{L,q}$ this grading (but not the name) has already been used in the proofs
of Theorems 2.3 and 2.4.

Returning to the cohomology vertex algebra $H_{d_{x}}(V)$ in the
case when $V$ is conformal and $x$ is an  eigenvector of $L_{0}$, we see that 
$$
   L_{-2}\b1\in Ker d_{x}\Rightarrow H_{d_{x}}(V)
\text{ is conformal,}
\eqno{(2.17)}
$$
$$
   L_{-2}\b1\in Im d_{x}\Rightarrow H_{d_{x}}(V)
\text{ is commutative.}
\eqno{(2.18)}
$$
Indeed, if $L_{-2}\b1\in Ker d_{X}$, then the operators $L_{i}\in End (V)$
descend to $H_{d_{X}}(V)$ due to (2.15).
 If, in addition, $L_{-2}\b1=  d_{x}(y)$,
then all $L_{i}$'s act on $H_{d_{x}}(V)$ trivially again due to (2.15). Hence 
$L_{0}$ acts on $H_{d_{X}}(V)$ trivially, each element of $H_{d_{X}}(V)$
is represented by a cocycle of conformal weight 0, and the $n$-th product
on $H_{d_{X}}(V)$ vanishes unless $n=-1$.

If $x\in V_{1}$, then $d_{x}(V_{n})\subset V_{n}$ for all $n$, and (2.18)
can be sharpened as follows:
$$
   L_{-2}\b1\in Im d_{x}\text{ and }x\in V_{-1} \Rightarrow  
 H_{d_{X}}(V)=H_{d_{X}}(V_{0}).
\eqno{(2.19)}
$$
 
In our previous work ([MSV], [MS1], [MS2]) we have dealt with conformal vertex
algebras having the following properties:   all conformal weights are
nonnegative; the conformal weight 0 component is a finitely
generated supercommutative ring and the corresponding multiplication coincides
with the restriction of the (-1)st multiplication. For example, 
$\Omega^{ch}_{X}$ is a sheaf of such vertex algebras
over a smooth manifold $X$: the conformal weight 0
component of $\Gamma(U,\Omega^{ch}_{X})$ is the algebra of differential forms
over $U\subset X$. In other words, there is a natural embedding 
$$
\Omega_{X}^{*}\iso\Omega^{ch}_{X,0}\subset \Omega^{ch}_{X},
\eqno{(2.20)}
$$
and it is this embedding that was invoked in Theorem 2.5B.

$H^{*}(X,\Omega^{ch}_{\BP^{N}})$ is also a vertex algebra of
this kind because its conformal weight 0 component equals  the cohomology
algebra $H^{*}(X)$. It is, therefore, natural to ask if there is a 
conformal vertex
algebra with nonnegative conformal weights so that the
(-1)-st multiplication identifies its conformal weight 0
component with the quantum cohomology of $X$. 

    The
quantum cohomology itself  is one such vertex algebra due to the equivalence
of categories reviewed above.
A more appealing possibility seems to be
 provided by $H_{D}(V_{L,q})$: it is a vertex algebra due to (2.16) 
because (1.1) is equivalent to
$$
D=\sum_{i=1}^{N}(\Psi^{i}_{-1}(e^{A^{i}}-e^{-\sum_{j}A^{j}}))_{(0)},
\eqno{(2.21)}
$$
and it is 
conformal because, as one easily checks, $D(L_{-2}e^{0})=0$.  

Even though Theorem 2.5A says that in this way we do not get anything new either, 
it
allows us to
 observe a curious phenomenon: $H_{D}(V_{L,q})$, $q\in\BC$,
 is a family of vertex
algebras over $\BC$ with fiber that equals $H^{*}(\BP^{N})$
 over any non-zero point and blows up to the  
non-commutative infinite dimensional vertex
algebra $H^{*}(\BP^{N},\Omega^{ch}_{\BP^{N}})$ over $0\in\BC$.

Rather unexpectedly, Theorem 2.5B turns out to be a by-product of 
the proof of Theorem 2.5A.

\bigskip 

{\bf 2.5.2} {\it Proof  Theorems 2.5A and B.}

By definition, the complex $(V_{L,q}, D)$ is
 the constant vector space
$V_{L}$ with differential $D$ polynomially depending on $q\in\BC$.
To make this more precise, observe that
 $V_{L}$ is graded by the function $ht$ defined in 2.1:
$$ 
V_{L}=\oplus_{n\geq 0}V_{L}^{n},
\eqno{(2.22)}
$$
where $V_{L}^{n}$ is a linear span of
$x\otimes e^{\sum_{i}m_{i}B^{i}+\sum_{i}n_{i}A^{i}}$ with
$ht(\sum_{i}n_{i}A^{i})=n$. The differential $D$ then breaks in a sum
$$
D=d_{+}+q^{N}d_{-},
\eqno{(2.23a)}
$$
so that
$$
d_{+}(V_{L}^{n})\subset V_{L}^{n+1},
\eqno{(2.23b)}
$$
$$
d_{-}(V_{L}^{n})\subset V_{L}^{n-N},
\eqno{(2.23c)}
$$
and
$$
(d_{+})^{2}=(d_{-})^{2}=[d_{+},d_{-}]=0.
\eqno{(2.23d)}
$$
Again by definition, the complex $(V_{L}, d_{+})$ coincides with
Borisov's complex $(V_{L}^{\Sigma},D)$. It follows from formulas (2.23a-d) 
and Theorem 1.3 that there
is a spectral sequence of the same type as (2.8), the 1st term and the
1st differential being as follows
$$
E_{1}=H^{*}(\BP^{N},\Omega^{ch}_{\BP^{N}})
\eqno{(2.24)}
$$
$$
d_{1}=q^{N}d_{-}:H^{*}(\BP^{N},\Omega^{ch}_{\BP^{N}})\rightarrow
H^{*}(\BP^{N},\Omega^{ch}_{\BP^{N}}),
$$
$$
d_{-}(H^{n}(\BP^{N},\Omega^{ch}_{\BP^{N}}))\subset
H^{n-N}(\BP^{N},\Omega^{ch}_{\BP^{N}}).
\eqno{(2.25)}
$$
Simply because $\text{dim} \BP^{N}=N$, the 2nd term equals
$$
\frac{H^{0}(\BP^{N},\Omega^{ch}_{\BP^{N}})}
{\text{Im}\{d_{-}:H^{N}(\BP^{N},\Omega^{ch}_{\BP^{N}})\rightarrow
H^{0}(\BP^{N},\Omega^{ch}_{\BP^{N}})\}}
$$
$$
\oplus 
\text{Ker}\{d_{-}:H^{N}(\BP^{N},\Omega^{ch}_{\BP^{N}})\rightarrow
H^{0}(\BP^{N},\Omega^{ch}_{\BP^{N}})\}\oplus\oplus_{i=1}^{N-1}
H^{i}(\BP^{N},\Omega^{ch}_{\BP^{N}}),
$$
and all higher differentials vanish. An argument similar to
(and simpler than) the one used in the proof of Lemma 2.4 shows that this
spectral sequence converges to $H_{D}(V_{L,q})$. Therefore
$$
H_{D}(V_{L,q})
$$
$$
=\frac{H^{0}(\BP^{N},\Omega^{ch}_{\BP^{N}})}
{\text{Im}\{d_{-}:H^{N}(\BP^{N},\Omega^{ch}_{\BP^{N}})\rightarrow
H^{0}(\BP^{N},\Omega^{ch}_{\BP^{N}})\}}
$$
$$
\oplus 
\text{Ker}\{d_{-}:H^{N}(\BP^{N},\Omega^{ch}_{\BP^{N}})\rightarrow
H^{0}(\BP^{N},\Omega^{ch}_{\BP^{N}})\}\oplus\oplus_{i=1}^{N-1}
H^{i}(\BP^{N},\Omega^{ch}_{\BP^{N}}).
\eqno{(2.26)}
$$

{\bf Lemma 2.6.} {\it There is $y\in H^{N}(\BP^{N},\Omega^{ch}_{\BP^{N}})$
such that $d_{-}(y)\in H^{0}(\BP^{N},\Omega^{ch}_{\BP^{N}})$
equals the Virasoro element $L_{-2}e^{0}$.}
\bigskip
This lemma allows us to complete the proof of  Theorems 2.5A and B
instantaneously. Our differentials come from elements of $V_{L}$ of conformal
weight 1, see (2.21); hence, due to Lemma 2.6, (2.18) and (2.19) apply:
 $H_{D}(V_{L,q})$
equals $H_{D}((V_{L,q})_{0})$, which is known (Theorem 2.4)
  to be equal to the
quantum cohomology. In particular, as follows from (2.26),  
$$
H^{i}(\BP^{N},\Omega^{ch}_{\BP^{N}})= 
H^{i}(\BP^{N},(\Omega^{ch}_{\BP^{N}})_{0}),\; 0<i<N,
$$
the latter space being canoncally isomorphic to
$H^{i}(\BP^{N},\Omega^{*}_{\BP^{N}})$ due to (2.20). Thus it remains to prove
Lemma 2.6.
\bigskip
{\it Proof of Lemma 2.6} 
To find an appropriate $y\in H^{N}(\BP^{N},\Omega^{ch}_{\BP^{N}})$ and calculate
$d_{-}(y)$ we need to take a  plunge in [MSV,B].

Let $x^{0}:x^{1}:\cdots :x^{N}$ be homogeneous coordinates on $\BP^{N}$ and 
 $b^{i}=x^{i}/x^{0}$. Consider the $N$-dimensional torus 
$\BT^{N}=\text{Spec}\BC[(b^{1})^{\pm 1},...,(b^{N})^{\pm 1}]\subset \BP^{N}$.

We shall need the following facts about the sheaf $\Omega^{ch}_{\BP^{N}}$.

First,
$$
\Gamma(\BT^{N}, \Omega^{ch}_{\BP^{N}})=
\BC[(b^{i}_{0})^{\pm 1},
b^{i}_{j-1},a^{i}_{j-1};\phi^{i}_{j},\psi^{i}_{j-1};1\leq i\leq N,j\leq 0],
\eqno{(2.27)}
$$
where $b^{i}_{j},a^{i}_{j-1}$ are even, $\phi^{i}_{j},\psi^{i}_{j-1}$ odd.

By letting $\text{deg}x^{i}_{j}=-j$, $x=b,a,\phi$ or $\psi$,
 we recover the grading by conformal weight.
By letting $\text{deg}b^{i}_{j}=\text{deg}a^{i}_{j}=0$,
$\text{deg}\phi^{i}_{j}=1$, $\text{deg}\psi^{i}_{j}=-1$ we get another grading, that
by {\it fermionic charge}. Therefore, $\Gamma(\BT^{N}, \Omega^{ch}_{\BP^{N}})$
is bigraded and this bigrading extends to the entire sheaf:
$$
\Omega^{ch}_{\BP^{N}}=\oplus_{m=-\infty}^{+\infty}\oplus_{n=0}^{+\infty}
\Omega^{ch,m}_{\BP^{N},n}.
\eqno{(2.28)}
$$
Next, we discuss ``tensor'' properties of $\Omega^{ch}_{\BP^{N}}$. We identify
$\Gamma(\BT^{N}, \Omega^{*}_{\BP^{N}})$ with 
$\Gamma(\BT^{N}, \Omega^{ch}_{\BP^{N}})_{0}$ by identifying $b^{i}$ with $b^{i}_{0}$
and $db^{i}$ with $\phi^{i}_{0}$. This identification extends to the
isomorphism (2.20).

The structure of higher conformal weight components is more complicated,
 but here is
what we can say about the component of conformal weight 1. Consider the following
  elements of $\Gamma(\BT^{N}, \Omega^{ch,0}_{\BP^{N},1})$:
$$
e_{ij}=b^{i-1}_{0}a^{j-1}_{-1}+\phi^{i-1}_{0}\psi^{j-1}_{-1},\; i,j\neq 1,
\eqno{(2.29a)}
$$
$$
e_{1j}=a^{j-1}_{-1},\; j\neq 1
\eqno{(2.29b)}
$$
$$
e_{i1}=-\sum_{l=1}^{N}b^{i-1}_{0}b^{l}_{0}a^{l}_{-1}
-\sum_{l=1}^{N}b^{i-1}_{0}\phi^{l}_{0}\psi^{l}_{-1}
$$
$$
-\sum_{l=1}^{N}b^{l}_{0}\phi^{i-1}_{0}\psi^{l}_{-1},\; i\neq 1.
\eqno{(2.29c)}
$$

It was checked in [MS1] III that these elements come from
$H^{0}(\BP^{N}, \Omega^{ch,0}_{\BP^{N},1})\subset
\Gamma(\BT^{N}, \Omega^{ch,0}_{\BP^{N},1})$ and that the Fourier components
of the corresponding fields span a Lie subalgebra
 of $\text{End}(\Omega^{ch}_{\BP^{N}})$
isomorphic to the  loop
 algebra $Lsl_{N+1}=sl_{N+1}\otimes\BC[t,t^{-1}]$.
 Therefore, 
$$
Lsl_{N+1}\hookrightarrow \text{End}(\Omega^{ch}_{\BP^{N}}),
\eqno{(2.30a)}
$$
so that
$$
sl_{N+1}\hookrightarrow H^{0}(\BP^{N}, \Omega^{ch,0}_{\BP^{N},1}),
E_{ij}\mapsto e_{ij}.
\eqno{(2.30b)}
$$
is a morphism of $sl_{N+1}$-modules, where $E_{ij}$, $i\neq j$,
$1\leq i,j\leq N+1$ are the standard generators of $sl_{N+1}$, and $sl_{N+1}$
operates on $H^{0}(\BP^{N}, \Omega^{ch,0}_{\BP^{N},1})$ by means of the 
composite map $sl_{N+1}\iso sl_{N+1}\otimes 1\subset \widehat{sl}_{N+1}$

Elements (2.29a-c) have fermionic charge 0. For the fermionic charge $N+1$
component there is an isomorphism:
$$
\Omega^{1}_{\BP^{N}}\otimes\Omega^{N}_{\BP^{N}}\iso
\Omega^{ch,N+1}_{\BP^{N},1}.
\eqno{(2.31)}
$$
  Over $\BT^{N}$ it is defined by the assignment 
$$
f_{i}(b^{1},...,b^{N})db^{i}\otimes(db^{1}\wedge db^{2}\wedge\cdots\wedge db^{N})
\mapsto 
f_{i}(b^{1}_{0},...,b^{N}_{0})
\phi^{i}_{-1}\phi^{1}_{0}\phi^{2}_{0}\cdots\phi^{N}_{0},
$$
$$
f_{i}(b^{1},...,b^{N})\subset\BC[(b^{1})^{\pm 1},...,(b^{N})^{\pm 1}].
$$
 Isomorphism (2.31) induces the isomorphism
$$
H^{N}(\BP^{N},\Omega^{1}_{\BP^{N}}\otimes\Omega^{N}_{\BP^{N}})\iso
H^{N}(\BP^{N},\Omega^{ch,N+1}_{\BP^{N},1}).
\eqno{(2.32)}
$$
By the Serre duality,
$$
H^{N}(\BP^{N},\Omega^{1}_{\BP^{N}}\otimes\Omega^{N}_{\BP^{N}})
\iso 
H^{0}(\BP^{N},\CT)^{*},
\eqno{(2.33)}
$$
where $\CT$ is the tangent sheaf. The Lie algebra
 $sl_{N+1}$ operates on $\BP^{N}$,
therefore there arises the map $sl_{N+1}\rightarrow H^{0}(\BP^{N},\CT)^{*}$,
 which
is well known to be  an isomorphism. Hence, (2.33) combined with (2.32)
 rewrites as follows
$$
 H^{N}(\BP^{N},\Omega^{ch,N+1}_{\BP^{N},1})
\iso 
sl_{N+1}.
\eqno{(2.34)}
$$
This map is an isomorphism of $sl_{N+1}$-modules, and it is not hard to find
 a Cech
cochain representing a highest weight vector of 
$H^{N}(\BP^{N},\Omega^{ch,N+1}_{\BP^{N}})_{1}$, that is, a non-zero vector $v$
 satisfying
$$
E_{ij}v=0,\; i<j.
\eqno{(2.35)}
$$
If we denote by $U_{i}$ the open subset of $\BP^{N}$ satisfying $x_{i}\neq 0$,
then $\{U_{0},...,U_{N}\}$ is an affine cover of $\BP^{N}$, so that
$\BT^{N}=U_{0}\cap U_{1}\cap...\cap U_{N}$. The $N$-th term
of the Cech complex equals, therefore, $\Gamma(\BT^{N},\Omega^{ch}_{\BP^{N}})$,
and it is an exercise to check that
$$
(b^{1}_{0})^{-1}(b^{2}_{0})^{-1}\cdots (b^{N-1}_{0})^{-1}(b^{N}_{0})^{-3}
\phi^{i}_{-1}\phi^{1}_{0}\phi^{2}_{0}\cdots\phi^{N}_{0}
\eqno{(2.36)}
$$
represents a highest weight vector of
 $H^{N}(\BP^{N},\Omega^{ch,N+1}_{\BP^{N},1})$.

Observe that another copy of $sl_{N+1}$ we have discovered earlier 
has $e_{1N+1}$ for its highest weight vector, see (2.29b,2.30). The  assertion 
crucial   for our proof is that $d_{-}$ sends one highest weight
 vector to another:
$$
d_{-}((b^{1}_{0})^{-1}b^{2}_{0})^{-1}\cdots b^{N-1}_{0})^{-1}b^{N}_{0})^{-3}
\phi^{i}_{-1}\phi^{1}_{0}\phi^{2}_{0}\cdots\phi^{N}_{0})=e_{1N+1}\in
H^{0}(\BP^{N}, \Omega^{ch,0}_{\BP^{N},1}).
\eqno{(2.37)}
$$
 Lemma 2.6 follows from (2.37) easily. To explain   this implication
we have to digress on   elementary  representation theory of
$Lsl_{N+1}$. 

\bigskip

Consider the  decomposition
$$
Lsl_{N+1}=L_{-}sl_{N+1}\oplus sl_{N+1}\oplus L_{+}sl_{N+1},
$$
where
$$
L_{\pm}sl_{N+1}=t^{\pm 1}\BC[t^{\pm 1}].
$$
Let $L_{\geq}sl_{N+1}=sl_{N+1}\oplus L_{+}sl_{N+1}.$ Any $sl_{N+1}$-module
becomes an $L_{\geq}sl_{N+1}$-module if the action of $sl_{N+1}$ is extended
to the entire $L_{\geq}sl_{N+1}$ by the requirement
$L_{+}sl_{N+1}\mapsto 0$. Therefore for any $sl_{N+1}$-module $U$ there arises
the {\it Weyl module}, denoted $\BW_{U}$ and defined as follows:
$$
\BW_{U}=\text{Ind}_{L_{\geq}sl_{N+1}}^{Lsl_{N+1}}U.
$$
The Weyl module induced from the trivial representation, $\BW_{\BC}$,
is well-known to be a conformal vertex algebra due to [FZ], see also
[K] 4.7. Therefore, it has vacuum
vector, $\b1$, and Virasoro element, $L^{aff}_{-2}\b1$. Other Weyl modules
are modules over $\BW_{\BC}$. This means, in particular, that Fourier
components $L^{aff}_{i}$ act on Weyl modules. The action of $L^{aff}_{0}$
is diagonalizable and defines a grading on each Weyl module also called the
grading by conformal weight. The aim of this digression was to formulate
the following well-known (and easily derived from the Kac-Kazhdan equations)
 assertion:
$$
I\subset \BW_{sl_{N+1}}\text{ is a {\it proper} $Lsl_{N+1}$-submodule }
\Rightarrow I\cap(\BW_{sl_{N+1}})_{2}=\{0\},
\eqno{(2.38)}
$$
where $\BW_{sl_{N+1}}$ stands for the Weyl module induced from the adjoint
representation, and $(\BW_{sl_{N+1}})_{2}$ is its conformal weight 2 component.

\bigskip

 Return to the proof of Lemma 2.6. Due to (2.30a), $H^{i}(\BP^{N},\Omega^{ch,m}_{\BP^{N}})$
is an $Lsl_{N+1}$-module for all $i$ and $m$. The component
$H^{N}(\BP^{N},\Omega^{ch,N+1}_{\BP^{N},1})$ is an $L_{\geq}sl_{N+1}$-module
isomorphic to $sl_{N+1}$, see (2.34), on which  $L_{+}sl_{N+1}$ acts trivially
because $H^{N}(\BP^{N},\Omega^{ch,N+1}_{\BP^{N},m})=0$ for all $m<1$. 
By the universality
property of induced modules, $\BW_{sl_{N+1}}$ maps onto the $Lsl_{N+1}$-submodule of
$H^{N}(\BP^{N},\Omega^{ch,N+1}_{\BP^{N}})$ generated by 
$H^{N}(\BP^{N},\Omega^{ch,N+1}_{\BP^{N},1})$.
Denote this submodule $\widehat{\BW}_{sl_{N+1}}$.

Similarly, $H^{0}(\BP^{N},\Omega^{ch,0}_{\BP^{N}})$ is an $Lsl_{N+1}$-module, and the
$Lsl_{N+1}$-submodule generated by $\b1$ is a quotient of $\BW_{\BC}$. This
quotient contains yet another submodule, the one generated by
$sl_{N+1}\hookrightarrow H^{0}(\BP^{N}, \Omega^{ch}_{\BP^{N}})_{1}$,
see (2.30b), to be denoted  $\widehat{\widehat{\BW}}_{sl_{N+1}}$.
  This submodule, again for the same reason, is
a quotient of $\BW_{sl_{N+1}}$. By definition, the above mentioned Virasoro
element $L^{aff}_{-2}\b1$ belongs to $(\widehat{\widehat{\BW}}_{sl_{N+1}})_{2}$

We are practically done. It is easy to derive from [B] that
$$
d_{-}:H^{N}(\BP^{N},\Omega^{ch,N+1}_{\BP^{N}})\rightarrow
H^{0}(\BP^{N},\Omega^{ch,0}_{\BP^{N}})
$$
is an $Lsl_{N+1}$-morphism. Equality (2.37) then means that 
$d_{-}(\widehat{\BW}_{sl_{N+1}})\subset \widehat{\widehat{\BW}}_{sl_{N+1}}$
 is  non-zero, and is therefore a quotient of 
$\BW_{sl_{N+1}}$ by a proper submodule. Due to (2.38)
$$
(d_{-}(\widehat{\BW}_{sl_{N+1}}))_{2}=(\BW_{sl_{N+1}})_{2}=
(\widehat{\widehat{\BW}}_{sl_{N+1}})_{2}.
$$
Hence $L^{aff}_{-2}\b1\in d_{-}(\widehat{\BW}_{sl_{N+1}})$. To complete the
proof of Lemma 2.6 it remains to check that the affine Virasoro
element,  $L^{aff}_{-2}\b1$, coincides with $L_{-2}\b1$ and this is easy.

\bigskip

Finally we have to prove (2.37). The difficulty with computation of
$$
d_{-}((b^{1}_{0})^{-1}(b^{2}_{0})^{-1}\cdots (b^{N-1}_{0})^{-1}(b^{N}_{0})^{-3}
\phi^{i}_{-1}\phi^{1}_{0}\phi^{2}_{0}\cdots\phi^{N}_{0})
$$
lies in that the operator $d_{-}$ is defined in terms of the vertex algebra $V_{L}$,
while 
$$
(b^{1}_{0})^{-1}(b^{2}_{0})^{-1}\cdots (b^{N-1}_{0})^{-1}(b^{N}_{0})^{-3}
\phi^{i}_{-1}\phi^{1}_{0}\phi^{2}_{0}\cdots\phi^{N}_{0}
$$
is an element of $\Gamma(\BT^{N},\Omega^{ch}_{\BP^{N}})$. The vertex algebra embedding
$$
\Gamma(\BT^{N},\Omega^{ch}_{\BP^{N}})\hookrightarrow V_{L},
$$
an important ingredient of Borisov's proof of Theorem 1.3, is determined by
the rules
$$
(b^{i}_{0})^{\pm}\mapsto e^{\pm B^{i}},
\phi^{i}_{0}\mapsto \Phi^{i}_{0}e^{B^{i}},
\psi^{i}_{-1}\mapsto \Psi^{i}_{-1}e^{-B^{i}},
\eqno{(2.39a)}
$$
$$
a^{i}_{-1}\mapsto A^{i}_{-1}e^{-B^{i}}-\Phi^{i}_{0}\Psi^{i}_{-1}e^{-B^{i}},
\eqno{(2.39b)}
$$
$$
x\mapsto X\;\Rightarrow L_{-1}x\mapsto L_{-1}X,
\eqno{(2.39c)}
$$
$$
x\mapsto X, y\mapsto Y\;\Rightarrow x_{(-1)}y\mapsto  X_{(-1)}Y.
\eqno{(2.39d)}
$$
These rules imply
$$
(b^{1}_{0})^{-1}(b^{2}_{0})^{-1}\cdots (b^{N-1}_{0})^{-1}(b^{N}_{0})^{-3}
\phi^{i}_{-1}\phi^{1}_{0}\phi^{2}_{0}\cdots\phi^{N}_{0}
\mapsto e^{-B^{N}}\Phi^{N}_{-1}\Phi^{1}_{0}\Phi^{2}_{0}\cdots \Phi^{N}_{0}.
$$

It follows from Borisov's proof of Theorem 1.3 that an element of $V_{L}$
representing the class of
$$
(b^{1}_{0})^{-1}(b^{2}_{0})^{-1}\cdots (b^{N-1}_{0})^{-1}(b^{N}_{0})^{-3}
\phi^{i}_{-1}\phi^{1}_{0}\phi^{2}_{0}\cdots\phi^{N}_{0}
$$
can be chosen to be equal to
$$
(\Psi^{N}_{-1}e^{A^{N}})_{(0)}
(\Psi^{N-1}_{-1}e^{A^{N-1}})_{(0)}\cdots
(\Psi^{1}_{-1}e^{A^{1}})_{(0)}
e^{-B^{N}}\Phi^{N}_{-1}\Phi^{1}_{0}\Phi^{2}_{0}\cdots \Phi^{N}_{0}.
$$
The formulas of 1.1-2 imply that  
$(\Psi^{i}_{-1}e^{A^{i}})_{(0)}$, $1\leq i\leq N-1$, simply erases $\Phi^{i}_{0}$.
Hence
$$
(\Psi^{N-1}_{-1}e^{A^{N}})_{(0)}
(\Psi^{N-2}_{-1}e^{A^{N-1}})_{(0)}\cdots
(\Psi^{1}_{-1}e^{A^{1}})_{(0)}
e^{-B^{N}}\Phi^{N}_{-1}\Phi^{1}_{0}\Phi^{2}_{0}\cdots \Phi^{N}_{0}
$$
$$
=\Phi^{N}_{-1}\Phi^{N}_{0}e^{-B^{N}+A^{1}+A^{2}+\cdots A^{N-1}}.
$$
The calculation of  the last operation is a little more tedious, but  also
 straightforward;
the result is this:
$$
(\Psi^{N}_{-1}e^{A^{N}})_{(0)}
(\Phi^{N}_{-1}\Phi^{N}_{0}e^{-B^{N}+A^{1}+A^{2}+\cdots A^{N-1}})
$$
$$
(\Psi^{N}_{-1}\Phi^{N}_{-1}\Phi^{N}_{0}-
\Phi^{N}_{-1}A^{N}_{-1}+
\frac{1}{2}\Phi^{N}_{0}A^{N}_{-2}+
\frac{1}{2}\Phi^{N}_{0}(A^{N}_{-1})^{2})e^{-B^{N}+A^{1}+A^{2}+\cdots A^{N}}.
\eqno{(2.40)}
$$
To complete our calculation we have to apply $d_{-}$ to this element. Observe
that this element comes from the interior of the cone spanned by 
$A^{1},...,A^{N}$ and has height $N$. It follows from the definition of the
spectral sequence and (1.1) or (2.21) that on this element $d_{-}$ equals
$$
-((\Psi^{1}_{-1}+\Psi^{2}_{-1}+\cdots+\Psi^{N}_{-1})
e^{-A^{1}-A^{2}-\cdots -A^{N}})_{(0)}.
$$
Indeed, it is precisely the component of Borisov's differential (1.1) that
decreases the height of the element (2.40). (By the way, it decreases it by $N$, which explains the assertion (2.23c).) Another calculation similar to
those performed  shows that
$$
-((\Psi^{1}_{-1}+\Psi^{2}_{-1}+\cdots+\Psi^{N}_{-1})
e^{-A^{1}-A^{2}-\cdots -A^{N}})_{(0)}.
$$
sends the element (2.40) to
$$
 A^{N}_{-1}e^{-B^{N}}-\Phi^{N}_{0}\Psi^{N}_{-1}e^{-B^{N}}.
$$
According to (2.39b), the latter element corresponds to $a^{N}_{-1}$ and hence to $e_{1N+1}$,
see (2.29b),
as desired. $\qed$

\bigskip\bigskip

\centerline{\bf \S 3. Deforming cohomology algebras of hypersurfaces
 in projective spaces}
 
\bigskip\bigskip

Let $\CL\rightarrow \BP^{N}$ be a degree $-n<0$ line bundle, 
$\CL^{\ast}\rightarrow \BP^{N}$ its dual,  $s\in\Gamma(\BP^{N},\CL^{\ast})$
a global section
so that its zero locus $Z(s)\subset \BP^{N}$ is a smooth   hypersurface. 
The way Borisov calculates the cohomology of the chiral de Rham complex over
$Z(s)$ is as follows.

Extend the lattice $(L, (.,.))$ introduced in 1.1 to the lattice
$(\hat{L}, (.,.))$ so that 
$$
\hat{L}=L\oplus\BZ A^{u}\oplus\BZ B^{u},\; (A^{u},B^{u})=1,(A^{u},L)=0,
(B^{u},L)=0.
$$

There arises the corresponding lattice vertex algebra $V_{\hat{L}}$.
Observe that any subset $L'\subset \hat{L}$ closed under addition gives
rise to the vertex subalgebra $V_{L'}\subset V_{L}$ 
 generated by $\gth_{\hat{L}}$ and
$Cl_{\hat{L}}$ from the highest weight vectors $e^{\beta},\beta\in L'$;
see 1.1-1.2.
 In our geometric
situation let $\hat{L}_{n}$ be the span of $B^{i}$ ($i=1,...,N$),
$B^{u}$ with arbitrary integral coefficients and $A^{i}$
($i=1,...,N$), $A^{u}$, $nA^{u}-A^{1}-\cdots - A^{N}$ with nonnegative
integral coefficients. 

The vertex  algebra $V_{\hat{L}_{n}}$ affords
a degeneration, $V_{\hat{L}_{n}}^{\Sigma}$, and includes in a family,
$V_{\hat{L}_{n},q}$, $q\in\BC$, in the same way the algebra $V_{L}$ did,
see 1.3, 2.1. To construct $V_{\hat{L}_{n}}^{\Sigma}$, consider 
 the following   $N+1$ elements of $\hat{L}$ :
$\xi_{1}=A^{1},\xi_{2}=A^{2},...,\xi_{N}=A^{N},
\xi_{N+1}=nA^{u}-A^{1}-A^{2}-\cdots
-A^{N}$. Define the   cone $\Delta_{i} $ to be
 the  set of all non-negative integral linear combinations of the elements
$\xi_{1},...,\xi_{i-1},\xi_{i+1},...,\xi_{N+1}, A^{u}$ and let
$\Sigma=\{\Delta_{1},...,\Delta_{N+1}\}$.  The vertex algebra
$V_{\hat{L}_{n}}^{\Sigma}$ is now defined by repeating word for word the
 definition of $V_{L}^{\Sigma}$ in 1.3. 

Similarly, the family
$V_{\hat{L}_{n},q}$, $q\neq 0$, is defined by repeating word for word the
 definition of $V_{L,q}$ in 2.1. This family extends ``analytically'' to
$q=0$ if $n\leq N+1$ and we again obtain an isomorphism 
$$
V_{L,0}=V_{\hat{L}_{n}}^{\Sigma}\text{ if }n<N+1.
\eqno{(3.1)}
$$
(The condition $n<N+1$ will be clarified below.)

Borisov's differential is as follows:
$$
D=\int\{\sum_{i=1}^{N}\Psi^{i}(z)(e^{A^{i}}-e^{nA^{u}-\sum_{j}A^{j}})(z)
+\Psi^{u}(z)(ne^{nA^{u}-\sum_{j}A^{j}}-e^{A^{u}})(z)\} .
\eqno{(3.2a)}
$$
(For the future use let us note that the right hand side of this equality
can be rewritten as a sum over lattice points:
$$
D=\int\{\sum_{i=1}^{N+1}\Psi^{\xi_{i}}(z)e^{\xi_{i}}(z) - 
\Psi^{u}(z)e^{A^{u}}(z)\},
\eqno{(3.2b)}
$$
where $\Psi^{\xi_{i}}=\Psi^{i}$ ($i\leq N$) and
$\Psi^{\xi_{N+1}}=n\Psi^{u}-\sum_{j}\Psi^{j}$.)

It is obvious that $D\in\text{End}(V_{\hat{L}_{n},q})$ and $D^{2}=0$;
 therefore there arise
the cohomology groups $H_{D}(V_{\hat{L}_{n},q} )$ and 
$ H_{D}(V_{\hat{L}_{n}}^{\Sigma})=H_{D}(V_{\hat{L}_{n},0} )$.
 
{\bf Theorem 3.1.} ([B])
$$ 
H_{D}(V_{\hat{L}_{n}}^{\Sigma}) =H^{*}(\CL,\Omega^{ch}_{\CL}).
$$
\bigskip

Borisov proposes to calculate the chiral de Rham complex over the
 hypersurface $Z(s)\subset \BP^{N}$
 by means of a certain Koszul-type resolution of the complex
$\Omega^{ch}_{\CL}$. The combinatorial data that determine 
$s\in\Gamma(\BP^{N},\CL^{\ast})$ 
consists of the finite set
$$
\Delta^{\ast}=\{\beta=B^{u}+\sum_{j=1}^{N}n_{j}B^{j}\text{ s.t. }
(\beta,\xi_{i})\geq 0,\; i=1,2,...,N+1\},
\eqno{(3.3)}
$$
and a function
$$
g: \Delta^{\ast}\rightarrow \BZ_{\geq}.
\eqno{(3.4)}
$$
Define
$$
K_{g}=\sum_{\beta\in \Delta^{\ast}}\int g(\beta)\Phi^{\beta}(z)e^{\beta}(z),
\eqno{(3.5)}
$$
where $\Phi^{\beta}=\Phi^{u}+\sum_{j}n_{j}\Phi^{j}$ provided
$\beta=B^{u}+\sum_{j}n_{j}B^{j}$. It is easy to see that
$$
K_{g}\in\text{End}(V_{\hat{L}_{n},q}),\; K_{g}^{2}=0, [K_{g},D]=0.
$$
 Therefore, there arise
the cohomology groups $H_{D+K_{g}}(V_{\hat{L}_{n},q} )$ and 
$ H_{D+K_{g}}(V_{\hat{L}_{n}}^{\Sigma})=H_{D+K_{g}}(V_{\hat{L}_{n},0} )$.

{\bf Theorem 3.2.} ([B])
$$ 
H_{D+K_{g}}(V_{\hat{L}_{n}}^{\Sigma}) =H^{*}(Z(s),\Omega^{ch}_{Z(s)}).
$$
\bigskip

All the vertex algebras in sight being topological (see the beginning of 2.2),
 Theorem 3.2 and the main result of [MSV]  give
$$
H^{*}(Z(s))=H_{Q_{0}}(H_{D+K_{g}}(V_{\hat{L}_{n}}^{\Sigma})),
\eqno{(3.6a)}
$$
or, equivalently,
$$
H^{*}(Z(s))=H_{D+K_{g}}(V_{\hat{L}_{n}}^{\Sigma})_{0},
\eqno{(3.6b)}
$$
where $H_{D+K_{g}}(V_{\hat{L}_{n}}^{\Sigma})_{0}$ stands for the kernel
of $L_{0}$.

This prompts the following  

{\bf Conjecture 3.3.} {\it If } $n<N+1$,
{\it  then the algebra } $H_{D+K_{g}}(V_{\hat{L}_{n},q} )_{0}$ 
{\it is isomorphic to the quantum cohomology algebra of} $Z(s)$.


Unfortunately we do not have a proof of this conjecture; we cannot
even prove that  $H_{D+K_{g}}(V_{\hat{L}_{n},q} )_{0}$ is  a 
deformation of $H^{\ast}(Z(s))$. What we know
is collected in the following 

{\bf Proposition 3.4.}
{\it (i) The element $e^{nA^{u}-\sum_{j}A^{j}}$ satisfies
$$
(D+K_{g})(e^{nA^{u}-\sum_{j}A^{j}})=0,
$$
and, therefore, determines an element of $H_{D+K_{g}}(V_{\hat{L}_{n},q} )_{0}$
for all $q$. If $q=0$, then this element, considered as an
element of $H^{\ast}(Z(s))$ (see (3.6b)), is proportional to  
 the cohomology class of a hyperplane section.

(ii) Due to (i), $e^{nA^{u}-\sum_{j}A^{j}}$ generates a subalgebra of
 $H_{D+K_{g}}(V_{\hat{L}_{n},q} )_{0}$ to be denoted $\CA_{q}$. This
subalgebra is a deformation of $\CA_{0}$.

(iii) If $Z(s)$ is a hyperplane (i.e. $n=1$), then Conjecture 3.3 is correct.

(iv) If $Z(s)$ is a non-degenerate quadric in $\BP^{3}$, then 
$H_{D+K_{g}}(V_{\hat{L}_{n},q} )_{0}$ is isomorphic
to $\BC[x,y]/(x^{2}-1,y^{2}-1).$. Hence Conjecture 3.3 is true in this case.}

 \bigskip

Since these results are by no means complete, we shall confine ourselves to
 sketching a proof of Proposition 3.4. The first part of assertion (i)
is a result of the obvious calculation using the formulas of 1.1-1.2.
The fact that at $q=0$ the element 
$e^{nA^{u}-\sum_{j}A^{j}}$ is proportional to
the cohomology class of a hyperplane section follows from Borisov's proof
of Theorem 3.2; this observation is completely analogous to the one made
in the end of 2.3. 

To prove (ii) observe that we have a constant family of vector spaces
$V_{\hat{L}_{n},q}$, $q\in\BC$, with differential $D+K_{g}$ depending
on $q$. At $q=0$ the complex $(V_{\hat{L}_{n},q},D+K_{g})$
degenerates in Borisov's complex $(V_{\hat{L}_{n}}^{\Sigma},D+K_{g})$.
As it always happens in situations of this kind, the differential
$d=D+K_{g}$ breaks in a sum $d=d_{-}(q)+d_{+}$ so that
$[d_{-}(q),d_{+}]=0$, $d_{-}(0)=0$, and $d_{+}$ equals Borisov's differential
on $V_{\hat{L}_{n}}^{\Sigma}$. There arises a spectral sequence converging
to $H_{D+K_{g}}(V_{\hat{L}_{n},q} )$ with  1-st term equal to
$H^{*}(Z(s),\Omega^{ch}_{Z(s)})$. The 2-nd term equals the cohomology
of the complex
$(H^{*}(Z(s),\Omega^{ch}_{Z(s)}), d^{(1)})$ with
$d^{(1)}=d_{-}(q))$. It remains to show that  
$$
\CA_{0}\subset \text{Ker} d^{(r)},\;
\CA_{0}\cap \text{Im} d^{(r)}=0,\; r\geq 1.
\eqno{(3.7)}
$$
All these spaces are subquotients of the subalgebra of $V_{\hat{L}_{n},0}$
generated by $e^{0}$, $e^{A^{i}}$ ($i=1,...,N$), $e^{A^{u}}$, and $\Phi^{i}_{0}$
($i=1,...,N$), $\Phi^{u}_{0}$, the product being equal to $_{(-1)}$. This
is a supercommutative algebra isomorphic to
$$
\BC[x_{1},...,x_{N},T,u;\Phi_{1},...,\Phi_{N},\Phi_{u}]/
(x_{1}x_{2}\cdots x_{N}T),
$$
where we let $x_{i}=e^{A^{i}}$, $u=e^{A^{u}}$, 
$T=e^{nA^{u}-\sum_{j}A^{j}}$, $\Phi_{i},\Phi_{u}$ being the corresponding
grassman variables. (All this is completely analogous to our discussion 
in the end of 2.1.)
Formula (3.2a) says that when
restricted to this space Borisov's
differential $D$ coincides with the Koszul differential associated with the
 sequence $x_{i}-T, u-nT$ ($i=1,...,N$) and our space quickly shrinks to
$$
\BC[T]/( T^{N+1}),
$$
on which (3.7) is obviously true at least when $r=1$. If $r\geq 2$, then the 
first part of (3.7) is obviously true and the second follows from a simple 
dimensional argument. 

Before turning to (iii)   
let us note that a quantum version of this argument gives:
$$
\CA_{q} \text{ is a quotient of } \BC[T]/(T^{N+1}-q^{N+1-n}n^{n}T^{n}).
\eqno{(3.8)}
$$
Indeed,  again by definition (as in the end of 2.1), 
the subalgebra of $V_{\hat{L}_{n},q}$
generated by $e^{A^{i}}$ ($i=1,...,N$), $e^{A^{u}}$, and $\Phi^{i}_{0}$
($i=1,...,N$), $\Phi^{u}_{0}$ is isomorphic to
 $$
\BC[x_{1},...,x_{N},T,u;\Phi_{1},...,\Phi_{N},\Phi_{u}]/
(x_{1}x_{2}\cdots x_{N}T-q^{N+1-n}u^{n}),
$$
and the restriction of $D$ to this supercommutative algebra coincides
with the Koszul differential associated with the regular sequence
$x_{i}-T, u-nT$ ($i=1,...,N$). The relation (3.8) follows at once. By the way,
the appearance of $N+1-n$ as a power of $q$ in (3.8) explains why the condition
$n<N+1$ was imposed in (3.1).

Return to the proof of (iii). In this case  the quantum cohomology
algebra is isomorphic to the algebra of functions on an $N$-point set.
Because of (ii), $\CA_{q}$ is isomorphic to $\BC[T]/p(T)$, $\text{deg}p(T)=N$,
and, because of (3.8), $p(T)$ divides $T^{N+1}-q^{N}T$. The latter has
no multiple roots. Hence $\CA_{q}$ is also the algebra of functions on an
$N$-point set. 

(iv) follows from the same spectral sequence that was used for
 the proof of (ii): due to (3.6b)
 $H^{\ast}(Z(s),\Omega^{ch}_{Z(s)})_{0}=\BC[x,y]/(x^{2},y^{2})$ and the 
elements $x,y$ are annihilated by all higher differentials because on the
 one hand
$x,y\in H^{1}(Z(s),\Omega^{ch}_{Z(s)})$ and on the other hand
it is true in general that all $d^{(r)}$, $r\geq 1$,
 send $H^{1}(Z(s),\Omega^{ch}_{Z(s)})$ to 0.
 The rest follows from (3.8), 
which in this case reads as follows:
$$T^{4}-4q^{2}T^{2}=0.
$$ $\qed$

\bigskip

{\bf Remarks.}  (i) By Corollary 9.3 of [G], the cohomology class $p$
 of a hyperplane section satsifies
in the quantum cohomology   of $Z(s)$  the relation
$$
p^{N}=q n^{n}p^{n-1}.
$$

  The amusing similarity between this equality and (3.8) suggests that
$\CA_{q}$ might be equal to $\BC[T]/(T^{N}-q^{N+1-n}n^{n}T^{n-1})$.

(ii)  Borisov's suggestion to treat the mirror symmetry as a flip
interchanging $A$'s and $B$'s seems to be working in our ``quantized''
situation as well. Compare (3.5) with (3.2b) to note
 that  $D$ and $K_{g}$ are sums over two
 sets of
lattice points defined by  self-dual condition (3.3).
 Hence the $A-B$ flip  changes $D$ to a similar differential to
be associated with the mirror partner of $Z(s)$ lying in another toric
manifold, see the next section.
 Of course the vertex algebra $V_{\hat{L}_{n}}$
bears a certain asymmetry, since not all elements of the type
$e^{\sum_{j}n_{j}A^{j}+n_{u}A^{u}}$  are allowed,
 but Borisov's ``transition to the
whole lattice'' (see Theorem 8.3 in [B]) and the above
spectral sequence seem to straighten things out.

\bigskip \bigskip

\newpage

\centerline{\bf \S 4. Quantum cohomology of toric varieties}

\bigskip\bigskip

Let us briefly explain how the constructions and results of section 2
carry over to an arbitrary smooth compact toric variety   of 
dimension $N$.  Each such variety 
is determined by a {\it complete regular}
fan in $L_{A}$. This and other relevant
concepts can be defined
as follows (see   [D, Bat] for details).

\bigskip

{\bf 4.1}
Let $I\subset L_{A}$. The cone generated by $I$ is said to be the
set of all non-negative integral combinations of elements of $I$
and is denoted $\Delta_{I}$.

A cone generated by part of a basis of $L_A$ is called {\it regular}.
 
A {\it complete regular fan} $\Sigma$ is defined to be a collection
of regular cones $\{\sigma_{1},...,\sigma_{s}\}$ so that the following
conditions hold:

(i) If $\sigma'$ is a face of $\sigma\in\Sigma$, then 
$\sigma'\in\Sigma$;

(ii) If $\sigma,\sigma'\in\Sigma$, then $\sigma\cap\sigma'$ is a face
of $\sigma$;

(iii) (the completeness condition)
 $L_{A}=\sigma_{1}\cup ...\cup\sigma_{s}$.

We   skip the construction of the smooth compact toric manifold
$X_{\Sigma}$ attached to a regular complete fan $\Sigma$ referring
the reader to [D], but  formulate Batyrev's result on
$H^{2}(X_{\Sigma}, \BR)$, see [Bat].

A function $\phi: L_{A}\rightarrow \BR$ is called 
 {\it piecewise linear} if its restriction to any cone in
$\Sigma$ is a morphism of abelian groups. Denote by $PL(\Sigma)$
the space of all piecewise linear functions.

Let $G(\Sigma)=\{\xi_{1},...,\xi_{n}\}$
 be the set of the generators of all 1-dimensional
cones in $\Sigma$. Since each piecewise linear function is determined
by its values on $\xi_{i}$ ($i=1,...,n$), $PL(\Sigma)$ is
 an $n$-dimensional real vector space. It contains the $N$-dimensional
subspace of globally linear functions; the latter is naturally isomorphic to 
$L_{B}\otimes_{\BZ}\BR$.

\bigskip

{\bf Theorem 4.1} ([Bat]) 
$$
H^{2}(X_{\Sigma}, \BR)=PL(\Sigma)/L_{B}\otimes_{\BZ}\BR.
$$

\bigskip

\bigskip

{\bf 4.2} Let us return to the vertex algebra $V_{L}$. Having fixed
 an arbitrary $\BR$-valued function $\phi$ on $L_{A}$, we proceed in 
much the same way as in 2.1. 

Define the linear automorphism 
$$
t_{\phi}: V_{L}\rightarrow V_{L} 
$$
by the formula
$$
t_{\phi}(x\otimes e^{\sum_{i}m_{i}B^{i}+\sum_{i}n_{i}A^{i}} )
=e^{-\phi(\sum_{i}n_{i}A^{i} )}
x\otimes e^{\sum_{i}m_{i}B^{i}+\sum_{i}n_{i}A^{i}}.
\eqno{(4.1)}
$$
Define $V_{L,\phi}$ to be the vertex algebra equal to $V_{L}$ as a vector
space with the following n-th product:
$$
(x\otimes e^{\sum_{i}m_{i}B^{i}+\sum_{i}n_{i}A^{i}})_{(n),\phi}
(y\otimes e^{\sum_{i}m'_{i}B^{i}+\sum_{i}n'_{i}A^{i}})
$$
$$=
t_{\phi}^{-1}(t_{\phi}(x\otimes e^{\sum_{i}m_{i}B^{i}+\sum_{i}n_{i}A^{i}} )
_{(n)}t_{\phi}(y\otimes e^{\sum_{i}m'_{i}B^{i}+\sum_{i}n'_{i}A^{i}} )).
\eqno{(4.2)}
$$
By definition,
$$
t_{\phi}: V_{L,\phi}\rightarrow V_{L}
$$
is a vertex algebra isomorphism. This provides us with a constant family
of vertex algebras parametrized by $\phi$ and we would like to study the
behavior of this family as $\phi$ tends to $\infty$. For this we have to impose
certain restrictions on $\phi$.

Following [Bat], call a piecewise linear function $\phi$  {\it convex} if
$$
\phi(x)+\phi(y)\geq\phi(x+y)\text{ all } x,y\in L_{A}.
\eqno{(4.3)}
$$
The cone of all convex piecewise linear functions descends to the cone
 in
$H^{2}(X_{\Sigma}, \BR)=PL(\Sigma)/L_{B}\otimes_{\BZ}\BR$, see Theorem 4.1. 
Denote this cone by $K(\Sigma)$ and its interior by
 $K^{0}(\Sigma)$.  $K^{0}(\Sigma)$ consists of classes of all {\it strictly
convex} piecewise linear functions, that is, of all
 those functions $\phi$
for which  equality in (4.1) is achieved if and only if $x$ and $y$
belong to the same cone in $\Sigma$.

We see immediately that

(i) if $\phi$ is convex piecewise linear, then the operations
$$
_{(n),\infty\phi}=\lim_{\tau\rightarrow +\infty}{_{(n),\tau\phi}},\;n\in\BZ
$$
are well defined and satisfy the Borcherds identities; denote the vertex
algebra arising in this way by $V_{L,\infty\phi}$;

(ii) if $\phi$ is strictly convex piecewise
 linear, then $V_{L,\infty\phi}$ is isomorphic
to Borisov's algebra $V_{L}^{\Sigma}$.

These assertions mean that the family $V_{L,\phi}$ produces a deformation of $V_{L}^{\Sigma}$ with base equal to the cone of
strictly convex piecewise linear functions. It is also immediate to see
that if $\phi-\phi'$ is a linear function, then the two deformations
$V_{L,\tau\phi}$  and $V_{L,\tau\phi'}$, $\tau\geq 0$, are equivalent.
Therefore we have obtained the family of vertex algebras 
$V_{L,\phi}$, $\phi\in K^{0}(\Sigma)$, which is a deformation of 
$V_{L}^{\Sigma}$ with base $K^{0}(\Sigma)$.

Denote by $QH^{\ast}_{\phi}(X_{\Sigma}, \BR)$ the quantum cohomology of
$X_{\Sigma}$ as defined in section 5 of [Bat].  Borisov's differential
is as follows
$$
D=\int\sum_{i=1}^{N}\Psi^{i}(z)(\sum_{j=1}^{n}(B^{i},\xi_{j})e^{\xi_{j}}(z)),
\eqno{(4.4)}
$$
where $\{\xi_{1},...,\xi_{n}\}$ is the set of generators of all 1-dimensional
cones in $\Sigma$.

\bigskip

{\bf Theorem 4.2}
$$
H_{Q_{0}+D}(V_{L,\phi})=QH^{\ast}_{\phi}(X_{\Sigma}, \BR).
$$

\bigskip

{\bf Sketch of Proof.} First of all, 
$$
(Q_{0})^{2}=0,\;D^{2}=0,\;[Q_{0},D]=0.
$$
(The first two of these assertions are obvious, the last one is obtained
in the same way as (2.4).)
Hence there arises a spectral sequence completely analogous to the one
used in 2.4. It converges and collapses:
$$
H_{Q_{0}+D}(V_{L,\phi})=H_{D}(H_{Q_{0}}(V_{L,\phi}));
$$
this is done in exactly the same way as in the proof of Theorem 2.4.

In part 1) of the proof
of Theorem 2.3  
the space $H_{Q_{0}}(V_{L,\phi})$ was shown to be equal to the group algebra $\BR[L_{A}]$ extended by grassman variables $\Phi^{i}_{0}$ ($i=1,...,N$).
Thus $H_{Q_{0}}(V_{L,\phi})$ is a Koszul complex, and 
the restriction of $D$ to this space equals the Koszul differential
associated with the sequence 
$$
\sum_{j=1}^{n}(B^{i},\xi_{j})e^{\xi_{j}},\; i=1,...,N.
\eqno{(4.5)}
$$
Therefore, $H_{D}(H_{Q_{0}}(V_{L,\phi}))$ is the corresponding
``Koszul cohomology''.

On the other hand,
 Batyrev defines $QH^{\ast}_{\phi}(X_{\Sigma}, \BR)$ to be the   polynomial ring $\BR[z_{1},...,z_{n}]$ modulo the sum of two ideals denoted $P(\Sigma)$
and $Q_{\phi}(\Sigma)$. It follows from the proof of either Theorem 9.5
or Theorem 8.4 in [Bat] that 
$$
\BR[L_{A}]= \BR[z_{1},...,z_{n}]/Q_{\phi}(\Sigma).
$$
Under this identification, the image of the ideal $P(\Sigma)$ in 
$\BR[L_{A}]$ is generated by the elements (4.5) as follows from the comparison
of (4.5) above and Definition 3.7 in [Bat]. The ring
$\BR[L_{A}]$ is Cohen-Macaulay; hence the sequence (4.5) is
regular. $\qed$ 

\bigskip 

\bigskip\bigskip


\centerline{\bf References}

\bigskip\bigskip

[Bat] V.~Batyrev, Quantum cohomology rings of toric manifolds, 
Journ\'ees de G\'eom\'etrie Alg\'ebrique d'Orsay (Orsay 1992), 
{\it Ast\'erisque}, {\bf 218} (1993), 9-34; alg-geom/9310004.   

[B] L.~Borisov, Vertex algebras and Mirror symmetry, 
math.AG/9809094. 

[D] V.I.~Danilov, The geometry of toric varieties, {\it Uspekhi Mat.Nauk}, 
{\bf 33}, No. 2 (1978), 85-134 (Russian); {\it Russian Math. Surveys}, 
{\bf 33}, No. 2 (1978), 97-154.

[FZ] I.B.~Frenkel, Y.~Zhu, Vertex operator algebras associated to
representations of affine and Virasoro algebra, {\it Duke Math. J.},
vol. {\bf 66}, No. 1, (1992), 123-168   

[G] A.~Givental, Equivariant Gromow-Witten invariants, {\it Internat. 
Math. Res. Notices}, {\bf 13} (1996), 613-663. 

[MSV] F.~Malikov, V.~Schechtman, A.~Vaintrob, Chiral de Rham complex, 
{\it Comm. Math. Phys.}, {\bf 204} (1999), 439-473.

[MS1] F.~Malikov, V.~Schechtman Chiral de Rham complex. II
{\it Amer. Math. Soc. Transl.} (2), vol. {\bf 194} (1999), 149-188

[MS2] F.~Malikov, V.~Schechtman Chiral Poincar\'e duality
{\it   Math. Res. Lett.}  vol. {\bf 6} (1999), 533-546

\bigskip\bigskip

F.M.:\ Department of Mathematics, University of Southern California, 
Los Angeles, CA 90089, USA; fmalikov\@mathj.usc.edu

V.S.:\ IHES, 35 Route de Chartres, 91440 Bures-sur-Yvette, France;  
vadik\@ihes.fr

\enddocument